\documentclass[psfig,a4paper,12pt,reqno]{amsart}
\usepackage{amsfonts,hyperref}
\usepackage{amssymb,amsmath,amsbsy,array,latexsym,color}

\theoremstyle{plain}
\begingroup
\newtheorem{theorem}    {Theorem}                       [section]
\newtheorem{lemma}      [theorem]        {Lemma}
\newtheorem{proposition}[theorem]        {Proposition}
\newtheorem{corollary}  [theorem]        {Corollary}
\theoremstyle{definition}
\newtheorem{example}    [theorem]        {Example}
\newtheorem{remark}     [theorem]        {Remark}
\newtheorem{definition} [theorem]        {Definition}
\theoremstyle{remark}

\numberwithin{equation}                                 {section}

\title{CONVEX FUNCTIONS ON SUB-RIEMANNIAN
MANIFOLDS. I}
\author{Kang-Hai Tan}
\address{Institute of Mathematics\\Fudan University, 200433\\ Shanghai, P. R. China}

\email{tankanghai2000@yahoo.com.cn,khtan@fudan.edu.cn}

\begin{document}

\begin{abstract}
 \small{We find a different approach to define convex functions in the sub-Riemannian setting.
 A function on a sub-Riemannian manifold is nonholonomically geodesic convex if its restriction to
 any nonholonomic (straightest) geodesic is convex. In the case of Carnot groups, this definition
 coincides with that by Danniell-Garofalo-Nieuhn
 (equivalent to that by Lu-Manfredi-Stroffolini). Nonholonomic
 geodesics are defined using the horizontal connection. A new
 distance corresponding to the horizontal connection has been
 introduced and near regular points proven to be equivalent to the
 Carnot-Carath\`{e}odory distance. Some basic properties of convex functions are
 studied. In particular we prove that any nonholonomically geodesic
 convex function locally bounded from above is locally Lipschitzian
 with respect to the Carnot-Carath\`{e}odory distance.
}
\end{abstract}

\maketitle

\section{Introduction}
Motivated by the role played by the theory of convex functions in
the theory of fully nonlinear partial differential equations-the
Monge-Amp\`{e}re equations for instance-in the Euclidean or
Riemannian case, mathematicians working in the field of subelliptic
partial differential equations have proposed several notions of
convexity of sets and functions in the setting of Carnot groups. The
notion of horizontal convexity (h-convexity for short) originally
formulated by Caffarelli, was rediscovered by Danielli, Garofalo and
Nhieu in \cite{DGN}. Roughly speaking, a subset $\Omega$ of a Carnot
group $G$ is said to be h-convex if the following condition holds:
if two points on an integral curve of some left invariant,
horizontal vector field on $G$ belong to $\Omega$, then the whole
segment of the integral curve between these two points is also
contained in $\Omega$. A function $f:\Omega\rightarrow\mathbb{R}$ is
h-convex if it is convex along the integral curves of the left
invariant, horizontal vector fields on $G$. The notion of horizontal
convexity in the sense of viscosity (v-convexity for short) was
proposed and studied by Lu, Manfredi and Stroffolli, see \cite{LMS}
for the Heisenberg group case and \cite{JLMS} for general Carnot
groups. Loosely speaking, an upper semicontinuous function
$f:\Omega\rightarrow\mathbb{R}$ defined on open subset of a Carnot
group $G$ is v-convex if the horizontal Hessian of test functions
touching $f$ from above is positive semidefinite. It turns out that
these two notions are equivalent, see
\cite{LMS,JLMS,BR,Ma,Wang,Ri,SY}. So far many fundamental properties
on horizontal convex functions have been obtained, see
\cite{DGN,LMS,DGNT,GT,GM,JLMS,BR,Ma,Wang,Ri,SY}. We in particular
remark that Rickly in \cite{Ri} proved that measurable h-convex
functions in Carnot groups are locally Lipschitzian with respect to
the Carnot-Carath\`{e}odory distance.

In this paper we will, from a more geometric viewpoint, define a
notion of convexity of functions on general sub-Riemannian
manifolds. This notion is based on the concept of nonholonomic
geodesics which researchers working in the field of nonholonomic
mechnics studied (see e.g. \cite{VF1,VF2,LD} and references
therein). Roughly speaking, nonholonomic geodesics in a
sub-Riemannian manifold are the straightest horizontal curves which
satisfy the equations of motion for the mechanical problem with a
quadratic Lagrangian and nonholonomic linear constraints which only
act by means of the reaction to them, i.e. in essence kinematically.
These curves can be characterized using a nonholonomic connection
(called horizontal connection in this paper). To be more precise,
let $(M, \Sigma,g_c)$ be a sub-Riemannian manifold where $\Sigma$ is
a subbundle of the tangent bundle $TM$ and $g_c$ a smooth inner
product on $\Sigma$. Without generality we assume that $g_c$ is the
restriction on $\Sigma$ of a Riemannian metric $g$ on $M$. For $X,
Y\in{\Gamma(\Sigma)}$, define $D_X Y := \mathcal{P}(\nabla_X Y )$
where $\nabla$ is the Levi-Civita connection of $g$ and
$\mathcal{P}$ denotes the projection onto $\Sigma$ with respect to
the metric $g$. Given a complement of $\Sigma$, i.e., a decompostion
$TM=\Sigma\oplus \Sigma^{\prime}$, $D$ depends only on $g_c$ (see
\cite{TY2}). Nonholonomic geodesics are those curves satisfying
$$
D_{\dot{\gamma}}\dot{\gamma}=0,\qquad
\dot{\gamma}\in{\Sigma_\gamma}.
$$
Nonholonomic geodesics are far different from sub-Riemannian
geodesics which are the shortest horizontal curves realizing the
Carnot-Carath\`{e}odory distance. In the literature the geometry
studying the truncated connection $D$ (resp. sub-Riemannian
geodesics) is called nonholonomic geometry, see e.g.
\cite{Sy,Va1,Va2}(resp. sub-Riemannian geometry, see e.g.
\cite{Mo}). However the connection $D$ plays an important role in
the study of sub-Riemannian geometry, see e.g. \cite{TY1,TY2} where
$D$ is used to define the horizontal mean curvature of hypersurfaces
of sub-Riemannian manifolds,  and \cite{Tan1} where we use $D$ to
define the notion of sublaplacian on sub-Riemannian manifolds. In
this paper we will use nonholonomic geodesics to give a notion of
convex functions in the setting of sub-Riemannian geometry. In
Section 3 we will use nonholonomic geodesics to define a new
distance $d_{\mathcal{H}}$. If a horizontal curve $\gamma$ consists
of smooth segements and moreover each segment is a nonholonomic
geodesic, we call $\gamma$ is a \textit{broken geodesic}. For
$p,q\in{M}$, let $d_{\mathcal{H}}(p,q)$ be the least length of all
broken geodesics connecting $p$ and $q$. Assume $\Sigma$ be bracket
generating and $M$ connected. In a neighborhood of each point we
will construct a horizontal frame such that the integral curves of
vector fields in this frame are nonholonomic geodesics. Thus we can
use the Chow connectivity theorem to prove for $p,q\in{M}$ there
exists at least a broken geodesic connecting $p$ and $q$, and
$d_{\mathcal{H}}$ is a distance. By definition $d_{\mathcal{H}}\geq
d_{c}$ where $d_c$ is the Carnot-Carath\`{e}odory distance. Indeed
by the ball-box theorem $d_c\leq d_{\mathcal{H}}\leq Cd_c$ locally
holds for a local constant $C$. For the first Heisenberg group
$\mathbb{H}^n$ we can prove $d_{\mathcal{H}}=d_c$. These results
tell us that the horizontal connection and its geodesics can provide
us another method to study sub-Riemannian geometry.

 To state our definition for convex functions in the setting of
 sub-Riemannian geometry, we first recall that in the Riemannian case a
function on a Riemannian manifold $M$ is convex if so is its
restriction to every Riemannian geodesic. Now let $f$ be a function
on the sub-Riemannian manifold $(M, \Sigma,g_c)$. We define that $f$
is nonholonomic geodesically convex (n-convex for short) if its
restriction to every nonholonomic geodesic is convex. Thus if
$\Sigma= TM$ this notion is just the Riemannian one. It is
remarkable that when $M$ is a Carnot group, this notion is just the
one by Danielli, Garofalo and Nhieu since now nonholonomic geodesics
are accidentally the integral curves of the left invariant,
horizontal vector fields. Monti and Rickly in \cite{MR} has proved
that in the Heisenberg group every function whose restriction to
every sub-Riemannian geodesic is convex must be constant. Thus in
our definition nonholonomic geodesics can not be replaced by
sub-Riemannian geodesics. We will also define a natural notion of
horizontal Hessian for smooth functions on sub-Riemannian manifolds
and then prove that a smooth function is n-convex  if and only if
its Hessian is positive semidefinite. It is natural to consider the
regularity problem for n-convex functions in general sub-Riemannian
manifolds. For Carnot groups, it is well known that each horizontal
convex function locally bounded from above is locally Lipschtizian.
The proof (see \cite{Ri}) of this fact depends only on a metric
property for Carnot groups $G$: there exist a constant $C$ and an
integer $N$ such that every two points $p,q\in{G}$ can be connected
by a broken geodesic, which is composed of $N$ nonholonomic
geodesics and each of them has length less than $Cd_c(p,q)$.
Fortunately this property locally holds for sub-Riemannian manifolds
$(M,\Sigma,g_c)$ where $\Sigma$ is regular and $M$ is connected.
Thus together with the local equivalence of $d_c$ and
$d_{\mathcal{H}}$ we prove local Lipschitz regularity for n-convex
functions with local bound from above.

One motivation to study convex functions on  such general
sub-Riemannian manifolds is the role played by the theory of convex
functions in the   study of structures of Riemannian manifolds. It
is well known that the existence of a convex function  on a
Riemannian manifold imposes strict limitations on its structure. We
expect that the theory of convex functions developed in this paper
could be used to study topology properties and structures of contact
manifolds, Riemannian submersions, and so on. This program will be
addressed in a forthcoming paper.

To end this introduction we give the structure of this paper. Some
basic facts about sub-Riemannian manifolds will be given in the next
section. We will adopt the viewpoint of Riemannian submersions to
study Carnot groups. The definition of horizontal Hessian on general
sub-Riemannian manifolds is new. Section 3 is devoted to
nonholonomic geodesics and the distance $d_{\mathcal{H}}$. In
Section 4 we study the n-convex notion and prove some basic
properties of n-convex functions. \vskip 10pt

{\small Acknowledgements. Theorem \ref{frame} is essentially due to
Dr. D. V. Tausk. The author thanks a lot for his helping in
understanding a similar construction in another context. The author
thanks Prof. Xiaoping Yang for his constant interest in this topic
and for many discussions. Finally he thanks Prof. Jiaxing Hong for
his encouragement and support when the author stayed in Fudan
University as a postdoc.}

\section{Horizontal connection and some associated differential
operators} Let $M$ be a smooth manifold of dimension $m$ endowed
with a smooth distribution (horizontal bundle) $\Sigma$ of dimension
$k$ with $k<m$. If we a prior equip $\Sigma$ with an inner product
$g_{c}$ (sub-Riemannian metric), we call $(M,\Sigma,g_{c})$ a
sub-Riemannian manifold with the sub-Riemannian structure
$(\Sigma,g_{c})$. If $\Sigma$ is integrable, it is just the
Riemannian geometry. We will assume $\Sigma$ is  not integrable. A
piecewise smooth curve $\gamma(t),t\in{[a,b]}$ in $M$ is horizontal
if $\dot{\gamma}(t)\in{\Sigma_{\gamma(t)}}$ a.e. $t\in{[a,b]}$. The
\textit{length} $\ell(\gamma)$ of the horizontal curve
$\gamma(t),t\in{[a,b]}$ is the integral
$\int_{a}^{b}g_{c}(\dot{\gamma}(t),\dot{\gamma}(t))dt$. Denote by
$\Sigma_{i}$ the set of all vector fields spanned by all commutators
of order $\leq i$ of vector fields in $\Sigma$  and let
$\Sigma_{i}(p)$ be the subspace of evaluations  at $p$ of all vector
fields in $\Sigma_{i}$.  We call $\Sigma$ satisfies the
Chow(H\"{o}rmander) condition or say $\Sigma$ is bracket generating
if for any $p\in{M}$, there exists an integer $l(p)$ such that
$\Sigma_{l(p)}(p)=T_{p}M$ (the least such $l$ is called the degree
of $\Sigma$ at $p$). If moreover $\Sigma_{i}$ is of constant
dimension (in a neighborhood of $p$) for all $i\leq l$, $\Sigma$ and
also $(M,\Sigma,g_c)$ are called regular (at $p$). If $M$ is
connected and $\Sigma$ satisfies the H\"{o}rmander condition, the
Chow connectivity theorem asserts that there exists at least one
piecewise smooth horizontal curve connecting two given points (see
\cite{Chow,Sub-Rie,Mo}), and thus $(\Sigma,g_{c})$ yields a metric
(called \textit{Carnot-Carath\'eodory distance}) $d_c$ by letting
$d_c(p,q)$ as the infimum among the lengths of all horizontal curves
joining $p$ to $q$. Those horizontal curves realizing $d_c$ are
sub-Riemannian geodesics. $(M,d_c)$ induces the same topology of
$M$. If $(M,d_c)$ is complete, then for any two points there exists
a shortest sub-Riemannian geodesic connecting them. Sub-Riemannian
geometry is anyway not  a trivial generalization of Riemannian
geometry. Some phenomena in sub-Riemannian geometry never appear in
Riemannian geometry. For example, in any neighborhood of a point $p$
there exists a point $q$ such that there are infinite many
sub-Riemannian geodesics connecting them, and  for some
sub-Riemannian manifolds there are sub-Riemannian geodesics (called
singular geodesics) not satisfying any ordinary differential
equation. We refer to \cite{Sub-Rie,Mo} for more about
sub-Riemannian geometry.

Sub-Riemannian manifolds are the common setting for control theory,
nonholonomic mechanics and many geometric structures such as CR
structures or contact metric structures, principal bundles with
connections, and Riemannian submersions. In this section we mainly
consider the horizontal connection which was originally introduced
to characterize equations of motion in nonholonomic Lagrangian
mechanics (\cite{Va1}).   We use the Einstein summation convention
for expressions with indices: if in any term the same index name
appears twice, as both an upper and a lower index, that term is
assumed to be summed over all possible values of that index (for $i,
j, r$ from 1 to $k$ and $a, b, c, d$ from 1 to $m$).
\begin{example}[Riemannian submersions] Let a smooth map $\pi$
between two Riemannian manifolds $\pi:(M, g)\rightarrow(B,
g^\prime)$ be a submersion, that is, $\pi_{\ast p}$ has maximal rank
at any point $p$ of $M$. Putting $\mathcal{V}_p =\ker(\pi_{\ast
p}$)for any $p\in M$, we obtain an integrable distribution
$\mathcal{V}$ which corresponds to the foliation of $M$ determined
by the fibres of $\pi$, since each $\mathcal{V}_p$ coincides with
the tangent space of $\pi^{-1}(x)$ at $p,$ $\pi(p)= x$.
$\mathcal{V}$ is called the vertical distribution whose section are
the so-called vertical vector fields. Let $\Sigma$ be the
complementary distribution of $\mathcal{V}$ determined by the
Riemannian metric $g$. Thus $(M, \Sigma,g|_\Sigma)$ is a
sub-Riemannian manifold. Given $X^\prime\in{\Gamma(TB)}$, the
horizontal vector field $X\in{\Gamma(\Sigma)}$ satisfying
$\pi_{\ast}(X) =X^\prime$ is called the horizontal lift of
$X^\prime$. Note that horizontal lifts of all vector fields of $B$
locally span $\Sigma$. We further assume that $\pi$ is a Riemannian
submersion, that is, moreover at each point of $p\in M$, $\pi_\ast$
preserves the length of the horizontal vectors. Since $p$ is a
submersion, $\pi_{\ast p}$ is a linear isomorphism between
$\Sigma_p$ and $T_{\pi(p)}B$ and $\pi_{\ast p}$ acts on $\Sigma_p$
as a linear isometry. For fundamental properties of Riemannian
submersions we refer to e.g. the book \cite{FLP}.
\end{example}
\begin{example}[Carnot groups]The most interesting models of
sub-Riemannian manifolds are Carnot groups (called also stratified
groups). The role played by Carnot groups in sub-Riemannian geometry
is the same as that by Euclidean spaces in Riemannian geometry (see
e.g. \cite{Mi,Be}). A Carnot group $G$ is a connected, simply
connected Lie group whose Lie algebra $\mathcal{G}$ admits the
grading $\mathcal{G}=V_1\oplus\cdots\oplus V_l$, with $[V_1,V_i]=
V_{i+1}$, for any $1\leq i\leq l-1$ and $[V_1,V_l] = 0$ (the integer
$l$ is called the step of $G$). Let $\{e_1,\cdots ,e_m\}$ be a basis
of $\mathcal{G}$ with $m=\sum_{i=1}^l\dim(V_i)$. Let
$X_i(g)=(L_g)_\ast e_i$ for $i=1,\cdots,m$ where $(L_g)_\ast$ is the
differential of the left translation $L_g(g^\prime)= gg^\prime$. We
call the system of left-invariant vector fields
$\Sigma:=V_1=\mathrm{span}\{X_1,\cdots,X_k\}$ the horizontal bundle
of $G$. If we equip $G$ an inner product $g$ such that
$\{X_1,\cdots,X_m\}$ is an orthonormal basis of $TG$, $(G,
\Sigma,g_c = g|_\Sigma)$ is a sub-Riemannian manifold. In $(G,
\Sigma,g_c)$ there exists a natural dilation homomorphism
$\delta_\lambda:\delta_{\lambda}p=\exp(\sum_{i=1}^{l}\lambda^{i}\xi_{i})$
for $p=\exp(\sum_{i=1}^{l}\xi_{i}), \xi_{i}\in{V_{i}}$. The most
simplest Carnot group is the Heisenberg group $\mathbb{H}^{n}$ which
is, by definition, simply $\mathbb{R}^{2n+1}$, with the
noncommutative group law
$$
pp^{\prime}=(x,y,t)(x^{\prime},y^\prime,t^{\prime})=(x+x^{\prime},y+y^{\prime},t+t^{\prime}+
\frac{1}{2}(\langle x^{\prime},y\rangle-\langle
x,y^{\prime}\rangle))
$$
 where we have let $x,
x^{\prime}, y, y^{\prime}\in{\mathbb{R}^{n}}, t,
t^{\prime}\in{\mathbb{R}}$. A simple computation shows that the
left-invariant vector fields $X_{j}(p)=\frac{\partial}{\partial
x_{j}}+\frac{1}{2}y_{j}\frac{\partial}{\partial t},
  X_{n+j}(p)=\frac{\partial}{\partial
y_{j}}-\frac{1}{2}x_{j}\frac{\partial}{\partial t}, j=1,\cdots,n, $
and $T=\frac{\partial}{\partial t}$ span the Lie algebra
($\mathbb{R}^{2n+1}$) of $\mathbb{H}^{n}$. Moreover
$[X_{j},X_{n+k}]=-T\delta_{jk}, j,k=1,\cdots,n,$ and all other
commutators are trivial. Note that the horizontal bundle
$\Delta=\textrm{span}\{X_{1},\cdots,X_{2n}\}$ is the kernel of the
1-form
$\eta=\frac{1}{2}dt+\frac{1}{2}\sum_{i=1}^{n}(x_{i}dy_{i}-y_{i}dx_{i})$
and the curvature form $\omega=d\eta=\sum_{i=1}^{n}dx_{i}\wedge
dy_{i}$ is the standard symplectic form in $\mathbb{R}^{2n}$. Thus a
piecewise smooth curve $\gamma(s)=(x(s),y(s),t(s)):[a,b]\rightarrow
\mathbb{H}^{n}$ is horizontal if and only if
\begin{equation}\label{horizontal}
\dot{t}(s)=\sum_{i=1}^{n}y_{i}(s)\dot{x_{i}}(s)-x_{i}(s)\dot{y}_{i}(s),\textrm{
when }\gamma\textrm{ is smooth at }s\in[a,b].
\end{equation}
\begin{lemma}\label{carnotsubmersion}Any Carnot group can be regarded as the source
manifold of a Riemannian submersion onto a standard Euclidean space
whose dimension equals to the dimension of the first layer of its
Lie algebra.
\end{lemma}
\proof Let $G$ be a Carnot group as above. Since the exponential map
is globally diffeomorphic, we usually identify $G$ with
$\mathbb{R}^n$ by the exponential map. It is easy to prove that
\begin{equation}\label{vector}
    X_j(x)=\frac{\partial}{\partial x_j}+\sum_{a=k+1}^mc_j^a(x)\frac{\partial}{\partial
    x_a},\qquad X_j(0)=e_j=\frac{\partial}{\partial
    x_j},j=1,\cdots,k
\end{equation}
where $c_j^a(x)=c_j^a(x_1,\cdots,x_k)$ are polynomials such that
$c_j^a(\delta_\lambda x)=\lambda^{w_a-1}c_j^a(x)$ where $w_a$ is the
weight of $x_a$ (\cite{Be,FS}), $a=k+1,\cdots,m$. For
$p=(x_1,\cdots,x_k,\cdots,x_m)\in G$, we define
$\pi(p)=(x_1,\cdots,x_k)\in{\mathbb{R}^k}$. Then by direct
computation it is easy to see that $\pi$ is a Riemannian submersion
from $(G, g)$ onto the Euclidean space $\mathbb{R}^k$ and the
horizontal lift of $\frac{\partial}{\partial x_j}$ is $X_j, j
=1,\cdots,k.$
\endproof
\noindent Although Lemma \ref{carnotsubmersion} is simple, via the
theory of Riemannian submersions it may simplify our discussions on
Carnot groups. For more information on Carnot groups we refer to
\cite{FS}.
\end{example}
Note that in a sub-Riemannian manifolds $(M,\Sigma,g_c)$ we always
can extend $g_c$ to a Riemannian metric $g$ in $M$ such that $TM$
can be $g$-orthogonally decomposed as $TM
=\Sigma\oplus\Sigma^\prime$ where $\Sigma^\prime$ is the
distribution complementary to $\Sigma$. We call such $g$ an
orthogonal extension of $g_c$. Obviously the orthogonal extension of
$g_c$ is not unique. We will use $\Gamma(\Sigma)$ to denote the set
of all smooth sections of $\Sigma$.
\begin{definition}[horizontal connection]Let $g$ be any orthogonal extension
 $g_c$ and let $\nabla$ be
the Levi-Civita connection with respect to $g$. We define the
horizontal connection $D$ on $\Sigma$ as
$$
D_XY=\mathcal{P}(\nabla_XY):=\sum_{i=1}^kg(\nabla_XY,X_i)X_i\qquad
\textrm{for any }X,Y\in{\Gamma(\Sigma)}
$$
where $\{X_1,\cdots,X_k\}$ is an orthonormal basis of $\Sigma$.
\end{definition}
It is obvious that $D$ is a linear connection on $\Sigma$. Since
$\Sigma$ is in general not integrable, such $D$ is called
nonholonomic connection in the literature.

\begin{proposition}[\cite{TY1}] Given a decomposition of $TM$,
$TM=\Sigma\oplus\Sigma^\bot$, $D$ is independent of the choice of
orthogonal extensions of $g_c$ such that $\Sigma$ is orthogonal to
$\Sigma^\bot$. Moreover $D$ is the unique nonholonomic connection on
$\Sigma$ satisfying
\begin{enumerate}
\item $Zg_c(X,Y)= g_c(D_ZX,Y)+ g_c(X,D_ZY)$;
\item $D_XY-D_YX =[X,Y]^\mathcal{H}$
\end{enumerate}
for any $X,Y,Z\in{\Gamma(\Sigma)}$, where
$[X,Y]^\mathcal{H}:=\mathcal{P}([X,Y])$ and $\mathcal{P}$ denotes
the (algebraic) projection onto $\Sigma$ with respect to the given
decomposition.
\end{proposition}
Thus $D$ depends only on the sub-Riemannian structure $(\Sigma,g_c)$
and the splitting of $TM$ and thus is ``almost" a sub-Riemannian
object. For many sub-Riemannian structures such as contact metric
structures, principal bundles with connections, and Riemannian
submersions the horizontal bundle $\Sigma$ has a canonical
complement. For these cases $D$ is canonical. Note that the
sub-Riemannian geometry of $(M,\Sigma,g_c)$, i.e., the geometry of
$(M,d_{c})$, depends only on the sub-Riemannian structure
$(\Sigma,g_c)$, not on complements of $\Sigma$ or extensions of
$g_c$. What should a canonical or good complement of $\Sigma$, and
then a canonical orthogonal extension of $g_c$ be? This may depend
on the context. In a try to find the sublaplacian for sub-Riemannian
manifolds the author in \cite{Tan1} proved the following statement
\begin{theorem}[\cite{Tan1}]Let $(M,\Sigma,g_c)$ be a sub-Riemannian manifold.
Then there exists a complement $\Sigma^\prime$ of $\Sigma$,
$TM=\Sigma\bigoplus\Sigma^{\prime}$, such that for this
decomposition there is an orthogonal extension $g$ of $g_c$ and an
orthonormal basis $\{T_1^\prime,\cdots,T^\prime_{m-k}\}$ of
$\Sigma^\prime$ satisfying
$$
\mathcal{P}(\nabla_{T^\prime_\beta}T^\prime_\beta)=0,\quad\beta=1,\cdots,m-k
$$
where $\nabla$ is the Riemannian connection of $g$. If $\Sigma$ is
strong-bracket generating(that is, for each $p\in{M}$ and each
nonzero horizontal vector $v\in{\Sigma_p}$ we have
$$
\Sigma_p+[V,\Sigma]_p=T_pM
$$
where $V$ is any horizontal extension of $v$), then such complement
is unique.
\end{theorem}
We will not use the last result in this paper. We cite it here just
to show that the selectivity of the splitting of the tangent bundle
and orthogonal extensions of $g_c$ is helpful in using $D$ to study
sub-Riemannian geometry. To simplify the discussion, unless
otherwise stated, we assume in the sequel that $g_c$ is the
restriction to $\Sigma$ of a given Riemannian metric $g$ on $M$ and
$TM$ admits the orthogonal decomposition, $TM
=\Sigma\oplus\Sigma^\bot$.

The horizontal connection $D$ induces the directional derivative of
a horizontal vector along a horizontal curve. Let $\{X_i\}_{i=1}^k$
be an orthonormal local basis of $\Gamma(\Sigma)$ and $Y$ be a
horizontal vector field along a horizontal curve $\gamma$. Then $Y$
can be locally written as $Y(t)= Y^j(t)X_j$ and the directional
derivative of $Y$ along $\gamma$ can be defined as
\begin{equation}\label{derivative}
\frac{D}{dt}Y(t)=\left(\dot{Y}^r(t)+Y^j(t)\dot{\gamma}^i(t)\Gamma^r_{ij}\right)X_r
\end{equation}
where $\Gamma^r_{ij}$ is Christoffel symbols such that $D_{X_i}X_j
=\Gamma_{ij}^rX_r$. Note that the horizontal connection is natural
in the sense that any isometry $\varphi$ between two sub-Riemannian
manifolds
$\varphi:(M^1,\Sigma^1,g_c^1)\rightarrow(M^2,\Sigma^2,g_c^2)$(i.e.
$\varphi$ is a diffeomorphism such that
$\varphi_\ast(\Sigma^1)=\Sigma^2$ and $\varphi^\ast g_c^2=g_c^1$)
takes the horizontal connection $D^1$ of $(M^1,\Sigma^1,g_c^1)$ to
the horizontal connection $D^2$ of
$(M^2,\Sigma^2,g_c^2):$
$$
\varphi_\ast(D^1_XY)=D^2_{\varphi_\ast(X)}(\varphi_\ast Y)
$$
for any $X,Y\in{\Gamma(\Sigma^1)}$.

\begin{lemma}\label{carnotconnection}In a Carnot group $G$, the horizontal connection has
the simple form
$$
D_UV=U(V^i)X_i\textrm{ for any }U,V=V^iX_i\in{\Gamma(\Sigma)}
$$
where $\{X_i\}_{i=1}^k$ is an orthonormal basis of the system of
left invariant, horizontal vector fields.
\end{lemma}
\proof Let $\{X_i\}_{i=1}^m$ be an orthonormal basis of left
invariant vector fields with respect to $g$. Then by the grading
condition of $\mathcal{G}$ we have
\begin{align*}
\Gamma_{ij}^r&=g_c(D_{X_i}X_i,X_r)
             =g(\nabla_{X_i}X_i,X_r)\\
             &=\frac{1}{2}\{g(X_i,[X_r,X_j])-g(X_j,[X_i,X_r])-g(X_r,[X_j,X_i])\}\\
             &=0
\end{align*}
Since $D_UV =U(V^i)X_i+U^jV^i\Gamma_{ij}^rX_r$, the statement
follows.
\endproof
\begin{definition}[horizontal gradient]For a smooth function $f:
M\rightarrow R$ we define its gradient as the horizontal vector
field $\nabla^\mathcal{H}f$ such that $g_c(X,\nabla^\mathcal{H}f)=
Xf$ holds for any horizontal vector field $X\in{\Gamma(\Sigma)}$.
\end{definition}
For every point $x\in M$ where $\nabla^\mathcal{H}f(x)=0$, by the
definition $\nabla^\mathcal{H}f(x)(\nabla^\mathcal{H}f(x))$ is the
horizontal direction along which $f$ increases (decreases) with the
fastest velocity.

\begin{definition}[horizontal divergence]For $X\in{\Gamma(\Sigma)}$, its
horizontal divergence is defined as
$\mathrm{div}^\mathcal{H}X=\sum_{i=1}^kg_c(D_{X_i}X,X_i)$ where
$\{X_i\}_{i=1}^k$ is an orthonormal local basis of $\Gamma(\Sigma)$.
It's clear that $\mathrm{div}^\mathcal{H}X$ is independent of the
choice of orthonormal bases.
\end{definition}
From Lemma \ref{carnotconnection}, it's easy to see that in a Carnot
group $G$, $\mathrm{div} X = \mathrm{div}^\mathcal{H}X$ for any
horizontal vector field $X$. Here div denotes the Riemannian
divergence with respect to $g$. We can give a notion of
sublaplacians on sub-Riemannian manifolds using the horizontal
gradient and divergence operators. We have given a detailed
discussion about this topic in \cite{Tan1}.
\begin{definition}[sublaplacian]The sublaplacian is defined as
$$
\Delta^\mathcal{H}:=\mathrm{div}^\mathcal{H}\circ\nabla^\mathcal{H}.
$$
\end{definition}
Next we define the horizontal Hessian of smooth functions on
sub-Riemannian manifolds.
\begin{definition}[horizontal Hessian]The horizontal Hessian of
a smooth function $f$ on $M$ is defined as
$$
\mathrm{Hess}^\mathcal{H}f(X,Y):=\frac{1}{2}\left\{g_c(DX(\nabla^\mathcal{H}f),Y)+
g_c(D_Y(\nabla^\mathcal{H}f),X)\right\}
$$
where $X,Y\in{\Gamma(\Sigma)}$.
\end{definition}
\begin{lemma}Let $\mathrm{Hess} f$ be the Riemannian Hessian of $f$. Then
\begin{equation}\label{hessianrep}
\mathrm{Hess}^\mathcal{H}f(X,Y)=\mathrm{Hess}
f(X,Y)-\frac{1}{2}\{B(X,Y)+B(Y,X)\}f
\end{equation}
for any $X,Y\in{\Gamma(\Sigma)}$, where $B(X,Y):=\nabla_XY-D_XY.$
\end{lemma}
It is obvious that $\mathrm{Hess}^\mathcal{H}f$ is a tensor on
$\Gamma(\Sigma)$. Note that the trace of the horizontal Hessian is
just the sublaplacian. The proof of the following statement is
trivial.

\begin{proposition}\label{carnothessian}Let $f$ be a smooth function on a Carnot group. Then
$$
\mathrm{Hess}^\mathcal{\mathcal{H}}f(X_i,X_j) =\frac{1}{2}(X_iX_jf +
X_jX_if),
$$
where $\{X_i\}_{i=1}^k$ as in Lemma \ref{carnotconnection}.
\end{proposition}
From Proposition \ref{carnothessian} we see that for Carnot groups
the horizontal Hessian is just the ordinary symmetrized horizontal
Hessian (e.g. \cite{DGN,LMS}. In particular, in Carnot groups the
sublaplacian is $\Delta^\mathcal{H}=\sum_{i=1}^k X_i^2$ where
$\{X_i\}_{i=1}^k$ as in Lemma \ref{carnotconnection}.

\section{Nonholonomic geodesics} In nonholonomic Lagrangian mechanics,
there are two well known approaches for the study of the constrained
mechanics: d'Alembertian nonholonomic mechanics and the variational
nonholonomic mechanics. The variational nonholonomic mechanics, also
called vakonomic mechanics by the Russian school, is to solve a
constrained variational problem: to find a curve $\gamma$ such that
$\gamma$ minimizes the Lagrangian functional $\int
L(\beta(t),\dot{\beta}(t))$ among all curves $\beta$ satisfying the
(nonholonomic) constraints $\dot{\beta}(t)\in{\Sigma_{\beta(t)}}$,
where $\Sigma$ is a (nonintegrable) distribution. When the
Lagrangian $L$ is regular, and quadratic with respect to the
velocity component, the required curve $\gamma$ is just a
sub-Riemannian geodesic, see e.g. \cite{Mo} for details. While the
dynamics studied by d'Alembertian nonholonomic mechanics is governed
by the Lagrange-d'Alembert principle. The principle states that the
equations of motion of a curve $q(t)$ in a configuration space are
obtained by setting to zero the variations in the integral of the
Lagrangian subject to variations lying in the constraint
distribution and that the velocity of the curve $q(t)$ itself
satisfies the constraints:
\begin{equation}\label{motioneq}
\dot{q}(t)\in{\Sigma_{q(t)}}\quad\textrm{ and}\quad -\delta
L:=\left(\frac{d}{dt}\frac{\partial L}{\partial
\dot{q}^i}-\frac{\partial L}{\partial q^i}\right)\delta q^i=0
\end{equation}
for all variations $\delta q$ such that $\delta q\in{\Sigma_q}$.
There is a huge literature on nonholonomic mechanics. We refer to
\cite{VF2,LD,Bl} and references therein.

The following theorem is well known, for its proof see e.g.
\cite{VF1,VF2}.
\begin{theorem}Assume $(M,\Sigma,g_c)$ be a sub-Riemannian manifold
where $g_c$ is the restriction to $\Sigma$ of a given Riemannian
metric $g$. If the Lagrangian $L$ is the kinetic energy
$L(q,\dot{q})=g(\dot{q},\dot{q})$, then by using the horizontal
connection $D$, the equations of motion \eqref{motioneq} can be
rewritten as
\begin{equation}\label{motioneq1}
D_{\dot{q}}\dot{q}=0.
\end{equation}
\end{theorem}
\begin{definition}[nonholonomic geodesic] On a sub-Riemannian
manifold $(M,\Sigma,g_c)$, any horizontal curve $\gamma$ satisfying
$D_{\dot{\gamma}}\dot{\gamma}=0$ is called a nonholonomic geodesic.
\end{definition}
The equation (3.2) of nonholonomic geodesics is a system of second
order differential equations formulated in local coordinates as
\begin{equation}\label{partialprayeq}
\frac{d^2q^c}{dt^2}+(\Gamma^\ast)_{ab}^c\frac{dq^a}{dt}\frac{dq^b}{dt}=0
\end{equation}
where $(\Gamma^\ast)_{ab}^c=\Gamma_{ab}^c+(\mu_i)_{a;b}(\mu_i)^c$,
$\Gamma_{ab}^c$ is the christoffel symbol,
$(\mu_i)_{a;b}=\frac{\partial(\mu_i)_a}{\partial
q^b}-\Gamma_{ab}^d(\mu_i)_d$ and $\mu_i$ are functions such that
$\Sigma$ is locally described by $\phi_i(q^a
,\dot{q}^a)=(\mu_i)_a(q)\dot{q}^a=0$, see \cite{Sy,LD}.

It is well known that Riemannian geodesics are the projections on M
of the integral curves of a spray. For nonholonomic geodesics, the
following theorem gives similar characterization.

\begin{theorem}[e.g.\cite{VF2,LD}]\label{partialspray} Nonholonomic geodesics are the
projections on $M$ of the integral curves of a vector field $\xi$ on
$\Sigma, \xi\in{\Gamma(T\Sigma)}$. Moreover there is an almost
product structure $(\mathbb{P}, \mathbb{Q})$ on $TM$ such that the
vector field $\xi$ can be explicitly given by projecting the
Riemannian spray $\xi^\prime$ to $T\Sigma$, that is,
$\xi(p)=\mathbb{P}(\xi^\prime(p))$ for $p\in\Sigma$.
\end{theorem}
We call the vector field $\xi$ in Theorem \ref{partialspray} a
\textit{partial spray}. Theorem \ref{partialspray} in particular
implies that for any point $x$ in $M$, and any vector
$v\in{\Sigma_x}$, there exists a unique smooth nonholonomic geodesic
$\gamma_v(t)$ such that $\gamma_v(0)=x$ and $\dot{\gamma_v}(0)=v$.
Note that any nonholonomic geodesic $\gamma$ has constant speed
because
$$
\frac{d}{dt}(g_c(\dot{\gamma},\dot{\gamma}))=2g_c(D_{\dot{\gamma}},\dot{\gamma})=0.
$$
From \eqref{partialprayeq} we have
\begin{lemma}
Given a point $x\in{M}$. For any $v\in{T_xM}$ and
$c,t\in\mathbb{R}$,
$$
\gamma_{cv}(t)=\gamma_v(ct)
$$
whenever either side is defined.
\end{lemma}
So for $x\in M$ there exists a neighborhood
$(0\in)\mathcal{E}_x\subset\Sigma_x$ such that for any
$v\in{\mathcal{E}_x}$, $\gamma_v(t)$ is defined on $[-2,2]$. We
define the \textit{horizontal( or partial) exponential}
$$
\exp_x^\mathcal{H}:\mathcal{E}_x\rightarrow M
$$
by $\exp_x^\mathcal{H}(v)=\gamma_v(1)$. Just as in the Riemannian
case, it is direct to prove that $\exp_x^\mathcal{H}$ is a
diffeomorphism from $\mathcal{E}_x$ onto a $k-$dimensional
submanifold containing $x$. Sometimes for $V=(x,v)\in{\Sigma}$, we
write $\exp^\mathcal{H}(V)$ (or $\gamma_V)$ for
$\exp_x^\mathcal{H}(v)$ (or $\gamma_v)$ when it is defined. Denote
by $\mathcal{E}\subset\Sigma$ the set on which a nonholonomic
geodesic $\gamma(t)$ is defined in $[-2,2]$. Then $\exp^\mathcal{H}$
is smooth on $\mathcal{E}$. The following theorem plays an important
role in this paper.
\begin{theorem}\label{frame}Given $y\in{(M,\Sigma,g_c)}$. There exists a neighborhood $\mathcal{O}\ni
y$ and a horizontal frame $\{E_1,\cdots,E_k\}$ on $\mathcal{O}$ with
$\|E_i\|^2:=g_c(E_i,E_i)=1$  such that any integral curve of $E_i$
is a nonholonomic geodesic in $\mathcal{O}$, $i=1,\cdots,k$.
\end{theorem}
\proof Let $\{X_1,\cdots,X_k\}$ be an orthonormal basis of $\Sigma$
in a neighborhood $\mathcal{U}$ of $y$ in $M$. For $X_i$ we choose a
hypersurface $S_i$ containing $y$ such that $X_i(x)$ is not in $T_x
S_i$ for all $x$ in $S_i\cap \mathcal{U}$. Now we can use the
horizontal exponential map $\exp^\mathcal{H}$ to obtain a ``local
coordinate system" near $S_i$ using nonholonomic geodesics with
initial velocities given by $X_i$. More explicitly, consider the
smooth map $\phi: (x,t)\rightarrow \exp^\mathcal{H}(tX_i(x))$, with
$x$ in $S_i$ and $t$ in an interval $(-\epsilon,\epsilon)$. By
taking $S_i$ and $\epsilon$ sufficiently small and using the inverse
function theorem one can prove that $\phi_i$ is a smooth
diffeomorphism onto an open neighborhood
$\mathcal{O}_i\subset\mathcal{U}$ of $y$ in $M$. Now we can extend
the vector field $X_i|_{S_i}$ to the open set $\mathcal{O}_i$ as
follows: given $p$ in $\mathcal{O}_i$ we pick $(x,t)$ with
$\phi_i(x,t)=p$ and define
$E_i(p)=\frac{d}{dt}\exp^\mathcal{H}(tX_i(x))$ in $T_pM$. Just by
definition it's clear  that  the integral curves of $E_i$ are
nonholonomic geodesics. Because $X_i$ has length 1, so is $E_i$. Set
$\mathcal{O}=\cap_{i=1}^k\mathcal{O}_i$. Taking $\mathcal{O}$
smaller if necessary, we get the desired frame $\{E_1,\cdots,E_k\}$
in $\mathcal{O}$.
\endproof
If a continuous horizontal curve $\gamma$ consists of segments each
of them a nonholonomic geodesic, we call $\gamma$ a \textit{broken
geodesic}. Thus broken geodesics are piecewise smooth. There is no
Hopf-Rinow type theorem for nonholonomic geodesics. Most couples of
points can not be connected by any nonholonomic geodesic. Through
the horizontal exponential it is easy to see that the set accessible
from a given point through nonholonomic geodesics is a smooth
submanifold of dimension $k$. But for broken geodesics we have
\begin{theorem}\label{connectivity}Let $(M,\Sigma,g_c)$ be a sub-Riemannian manifold.
If  $\Sigma$ satisfies the Chow condition and $M$ is connected, for
any two points $p,q\in{M}$ there exists at least a broken geodesic
joining $p$ to $q$.
\end{theorem}
\proof For any $p\in M$, let $\mathcal{A}_p$ be the set of points
which are reachable via broken geodesics starting from $p$. To prove
the statement it is sufficient to prove the openness of
$\mathcal{A}_p$ since $M$ is connected and $\mathcal{A}_p$ is
trivially not empty. For $y\in{\mathcal{A}_p}$, by Theorem
\ref{frame} in a neighborhood $\mathcal{O}$ of $y$ there exists a
frame $\{E_1,\cdots,E_k\}$ such that integral curves of $E_i$ are
nonholonomic geodesics. Since
$\Sigma|_\mathcal{O}=\mathrm{span}\{E_1,\cdots,E_k\}$ is bracket
generating, by Chow connectivity theorem there exists a neighborhood
of $y$ such that any point in this neighborhood can be connected by
a piecewise smooth horizontal curve starting from $y$ and each piece
is an integral curve of $E_i$.
\endproof
Now we define a new distance through broken geodesics.
\begin{definition}Let $(M,\Sigma,g_c)$ be a sub-Riemannian manifolds
with $\Sigma$ bracket generating and $M$ connected. For $p,q\in{M}$,
the distance $d_{\mathcal{H}}(p,q)$ is defined as the least length
among all broken geodesics connecting $p$ and $q$.
\end{definition}
From Theorem \ref{connectivity} we know $d_{\mathcal{H}}$ is well
defined and by definition $d_c\leq d_{\mathcal{H}}$. To get more
about the distance $d_{\mathcal{H}}$ we introduce some notations.
Let $Y_1$ and $Y_2$ be smooth vector fields, with local flows
$\Psi_i(t)=\exp(tY_i)$. Then for small $t$,
$$
\Psi_1(-t)\circ\Psi_2(-t)\circ\Psi_1(t)\circ\Psi_2(t)(p)=p+t^2[E_1,E_2](p)+O(t^2).
$$
Write $[Y_1(t),Y_2(t)]$ for
$\Psi_1(-t)\circ\Psi_2(-t)\circ\Psi_1(t)\circ\Psi_2(t)$. For $y\in
M$, let $\{E_1,\cdots,E_k\}$ in $\mathcal{O}(\ni y)$ be the frame of
norm 1 in Theorem \ref{frame}. For multi-indices
$I=(i_1,\cdots,i_r), 1\leq i_j\leq k$, define vector fields $E_I$
inductively by $E_I=[E_{i_1},E_J],$ where $J=(i_2,\cdots,i_r)$. We
write $i_1J=I$ and denote the length of a multi-index $I$ by $|I|$,
so $|J|=r-1$. Similarly define flows
$\Psi_I(t)=[\Psi_{i_1}(t),\Psi_J(t)]$ as above for $Y_1,Y_2$. Note
that $\Psi_I(t)=\textbf{1}+t^rX_I+O(t^{r+1})$ and $\Psi_I(t)y$ is a
concatenation of $(3\cdot2^{r-1}-1)$ nonholonomic geodesics, each
one of length $\epsilon$ if $|t|\leq\epsilon$. If $\Sigma$ is
bracket generating, we can select a  frame for the entire tangent
bundle (of $\mathcal{O}$) among the $E_I$. We choose such a frame
and relabel it
$\{E_1,\cdots,E_{n_1},E_{n_1+1},\cdots,E_{n_2},E_{n_2+1},\cdots,E_m\}$
where $n_1=k$, $\{E_1,\cdots,E_{n_i}\}$ spans
$\Sigma_i:=\Sigma_{i-1}+[\Sigma_{i-1},\Sigma],\Sigma_1=\Sigma,
i=1,\cdots,l$ and $n_l=m$. $(k,n_1,\cdots,n_l)$ is the growth vector
of $\Sigma$ at $y$. For each chosen $E_i$ of the form $E_I$, let
$w_i$ be the length $|I|$. We also relabel flows $\Psi_I$ as
$\Psi_i$, $i=1,\cdots,m$. Coordinates $x_1,\cdots,x_m$ are said to
be linearly adapted to $\Sigma$ in $\mathcal{O}$ if $\Sigma_i$ is
annihilated by the differentials $dx_{n_i+1},\cdots,dx_m$ in
$\mathcal{O}$. The \textit{weighted box} of size $\epsilon$ is the
set
$$
\mathrm{Box}^w(\epsilon)=\{x\in{\mathbb{R}^m}:|x_i|\leq\epsilon^{w_i},i=1,\cdots,m\}.
$$
Define the map
$F^y(t_1,\cdots,t_m)=\Psi_m(t_m)\circ\cdots\circ\Psi_1(t_1)(y):\mathbb{R}^m\rightarrow\mathcal{O}$.

In sub-Riemannian geometry, the ball-box theorem is well-known. The
following theorem is essentially proven by \cite[p.27-34]{Mo}, see
also \cite{Be}.
\begin{theorem}\label{ballbox}Let $(M,\Sigma,g_c)$ be a sub-Riemannian manifold (with
$\Sigma$ bracket generating and $M$ connected). Then for a regular
point $y_0\in M$ there exist a neighborhood $\mathcal{O}$ and
linearly adapted coordinates $x_1,\cdots,x_m$ such that for any
$y\in{\mathcal{O}}$ there exist positive continuous functions
$c(y)<C(y)$, $\epsilon_0(y)$ and an integer
$N=\sum_{r=1}^l(n_r-n_{r-1})(3\cdot2^{r-1}-1)$ such that for all
$\epsilon<\epsilon_0$,
\begin{equation}\label{ballbox1}
    B(y,\frac{c}{C}\epsilon)\subset F^y(\mathrm{Box}^s(N^{-1}\epsilon))\subset
    B(y,\epsilon)
\end{equation}
where
$\mathrm{Box}^s(\epsilon):=\{(t_1,\cdots,t_m)\in{\mathbb{R}^m}:|t_i|\leq\epsilon\}$
is the standard $\epsilon-$cube and
$B(y,\epsilon):=\{p\in{\mathcal{O}}:d_c(y,p)\leq\epsilon\}$ denotes
the Carnot-Carath\`{e}odory ball centered at $y$.
\end{theorem}
\proof Note that because $\Sigma$ is regular at $y_0$, $N$ is a
constant in a neighborhood $\mathcal{O}$ of $y_0$. From the proof of
the ball-box theorem given in \cite[p.27-34]{Mo}, we know that near
a regular point there exist linearly adapted coordinates
$x_1,\cdots,x_m$ and positive constants $c<C$, $\epsilon_0$ and an
integer $N=\sum_{r=1}^l(n_r-n_{r-1})(3\cdot2^{r-1}-1)$ such that for
all $\epsilon<\epsilon_0$,
\begin{equation}\label{ballbox2}
\mathrm{Box}^w(c\epsilon)\subset
F^y(\mathrm{Box}^s(N^{-1}\epsilon))\subset B(y,\epsilon)\subset
\mathrm{Box}^w(C\epsilon)
\end{equation}
and $c,C$ and $\epsilon_0$ continuously depend on
$y\in{\mathcal{O}}$. We first use \eqref{ballbox2} replacing
$\epsilon$ by $\frac{c}{C}\epsilon$ to get
$B(y,\frac{c}{C}\epsilon)\subset \mathrm{Box}^w(c\epsilon)$. Using
\eqref{ballbox2} again we obtain \eqref{ballbox1}.
\endproof
\begin{remark}If $y_0$ is not a regular point, the constants $c,C$
and $\epsilon_0$ are in general not continuous functions of $y$.
\end{remark}

 The ball-box theorem implies the Chow connectivity theorem
and that the topology of $(M,d_c)$ is the same as the original one.

\begin{corollary}\label{corcontrol}Let $(M,\Sigma,g_c)$ be a sub-Riemannian manifold
(with $\Sigma$ bracket generating and $M$ connected). Then for any
regular  point $y_0\in M$ there exists a sub-Riemannian ball
$B(y_0,\bar{\epsilon}_0)$ centered at $y_0$ with radius
$\bar{\epsilon}_0$, and two constants $C_1,C_2$ such that for any
two points $p,q\in{B(y_0,\bar{\epsilon}_0)}$
\begin{enumerate}
  \item\label{control}
 there exists a
broken geodesic connecting $p$ and $q$, which consists of $N$
nonholonomic geodesics, each one of length less than $C_1d_c(p,q)$;
  \item $d_c(p,q)\leq d_{\mathcal{H}}(p,q)\leq C_2d_c(p,q)$.
\end{enumerate}
\end{corollary}
\proof Let $\mathcal{O}^\prime$ be a neighborhood of $y_0$ such that
$\bar{\mathcal{O}^\prime}\subset\mathcal{O}$ where $\mathcal{O}$ is
the neighborhood in Theorem \ref{ballbox}. In
$\bar{\mathcal{O}}^\prime$ set $c_1$ be the maximum of $C(y)$, $c_2$
the minimum of $c(y)$ and $\epsilon_0^\prime$  the minimum of
$\epsilon_0(y)$. Take a sub-Riemannian ball
$B(y_0,\bar{\epsilon}_0)\subset\mathcal{O}^\prime$ with
$\bar{\epsilon}_0<\frac{c_2}{4c_1}\epsilon_0^\prime$. Then by
Theorem \ref{ballbox} for any $y\in{B(y_0,\bar{\epsilon}_0)}$ and
$\epsilon<\epsilon_0^\prime$ we have
\begin{equation}\label{uniform}
\overline{B(y,\frac{c_2}{2c_1}\epsilon)}\subset
B(y,\frac{c_2}{c_1}\epsilon)\subset
F^y(\mathrm{Box}^s(N^{-1}\epsilon)).
\end{equation}
For any $p,q\in{B(y_0,\bar{\epsilon}_0)}$, because
$d_c(p,q)=\frac{c_2}{2c_1}\bar{\epsilon}$ for
$\bar{\epsilon}:=\frac{2c_1d_c(p,q)}{c_2}<\frac{4c_1\bar{\epsilon}_0}{c_2}<\epsilon_0^\prime$,
 by \eqref{uniform} we obtain
$q\in{F^p(\mathrm{Box}^s(N^{-1}\bar{\epsilon}))}$, that is, $q$ is
the endpoint of  a broken geodesic $\gamma$ starting from $p$ and
consisting of $N$ nonholonomic geodesics, each one of length less
than $N^{-1}\bar{\epsilon}=C_1d_c(p,q)$ with $C_1=\frac{2c_1}{Nc_2}$
(Take $\bar{\epsilon}_0$ smaller if necessary such that for any
$p\in{B(y_0,\bar{\epsilon}_0)}$, $F^p$ is well defined in
$\mathcal{O}^\prime$). We also get $d_{\mathcal{H}}(p,q)\leq
C_2d_c(p,q)$ with $C_2=\frac{2c_1}{c_2}$.
\endproof
\begin{remark}\quad
\begin{enumerate}
                \item The condition ``with $\Sigma$ bracket generating and $M$
connected"  in Theorem \ref{ballbox} and Corollary \ref{corcontrol}
is not necessary, since we assume $y_0$ is regular, $d_c$ and
$d_{\mathcal{H}}$ can be defined in $\mathcal{O}$.
                 \item From the proof of Corollary \ref{corcontrol} we see that
if moreover $M$ is compact then the constants $C_1,C_2$ and
$\bar{\epsilon}_0$ can be taken universal constants.
                \item If $\Sigma$ is regular and $M$ connected, the
topology of $(M,d_{\mathcal{H}})$ is the same as that of $(M,d_c)$.
                \item As pointed out in Section 2, the horizontal
connection depends on the choice of the complement of $\Sigma$. So
different choices of splitting lead to different sets of
nonholonomic geodesics. However Corollary \ref{corcontrol} holds for
all choices with $C_1$ and $C_2$ possibly changed up to the choice.
              \end{enumerate}
\end{remark}
In the remainder of this section we will consider nonholonomic
geodesics and $d_{\mathcal{H}}$ for some special sub-Riemannian
manifolds, in particular for Carnot groups.
\begin{theorem}\label{thmrieman}Let $(M,\Sigma,g_c)$ be a sub-Riemannian manifold where
$g_c = g|_\Sigma$ for a Riemannian metric $g$ on $M$. If for any
horizontal vector field $X\in{\Gamma(\Sigma)}$ and any vertical
vector field in the orthogonal complement,
$Y\in{\Gamma(\Sigma^\bot)}$ one has $[X,Y]\in{\Gamma(\Sigma^\bot)}$,
then a horizontal curve is a nonholonomic geodesic of
$(M,\Sigma,g_c)$ if and only if it is a Riemannian geodesic of $(M,
g)$.
\end{theorem}
 \proof Let $\gamma$ be a
nonholonomic geodesic (assume not a constant). Since $\gamma$ has
constant speed it is regular and one can extend $\dot{\gamma}$ to a
smooth horizontal vector field $Y$ in a small neighborhood
$\mathcal{U}$ of $\gamma(t_0)$ such that
$\|Y\|^2:=g_c(Y,Y)\geq\frac{1}{2}c$ in $\mathcal{U}$, where
$c:=g_c(\dot{\gamma},\dot{\gamma})$. Then $X:=c\frac{Y}{\|Y\|}$ is
also a smooth extension of $\dot{\gamma}$ with constant norm $c$. We
claim that the projection on $\Sigma^\bot$ of $\nabla_XX$ vanishes,
i.e., $(\nabla_XX)^\bot=0$ in $\mathcal{U}$. In fact, for any
vertical vector field $W\in{\Gamma(\Sigma^\bot)}$, we have
\begin{align*}
 g((\nabla_XX)^\bot,W)&=-g(\nabla_XW,X)\\
                      &=-g(\nabla_WX,X)+ g([W, X],X)\\
                      &=-Wg(X,X)\\
                      &=0
\end{align*}
 where we
used the condition that $[W, X]\in{\Gamma(\Sigma^\bot)}$ and $g(X,
X)\equiv c$. Thus $(\nabla_{\dot{\gamma}}\dot{\gamma}(t))^\bot=0$ in
a small neighborhood of $t_0$. Thus $\dot{\gamma}$ is a Riemannian
geodesic. The inverse follows from a similar argument.
\endproof

Because for a Riemannian submersion the Lie bracket of a horizontal
(projectable) vector field with a vertical vector field is vertical
(see e.g. [6]), we have the following statement.

\begin{corollary}\label{corosubmersion}Let p be a Riemannian submersion $\pi:(M,g)\rightarrow(B,
g^\prime)$ and let $\Sigma$ be the horizontal bundle of $(M, g)$.
Then a horizontal curve is a nonholonomic geodesic of
$(M,\Sigma,g_c= g|_\Sigma)$ if and only if it is a Riemannian
geodesic of $(M,g)$. Thus the set of nonholonomic geodesics of
$(M,\Sigma,g_c)$ consists of horizontal lifts of all Riemannian
geodesics of $(B, g^\prime)$ because the horizontal lift of every
Riemannian geodesic of $(B, g^\prime)$ is a Riemannian geodesic of
$(M,g)$.
\end{corollary}
\begin{corollary}\label{corocarnotgeodesic} A horizontal curve in a Carnot group is a
nonholonomic geodesic if and only if it is an integral curve of a
left invariant horizontal vector field.
\end{corollary}
\proof Let $G$ be a Carnot group. The submersion $\pi$ given in the
proof of Lemma \ref{carnotsubmersion} satisfies
$\pi(p.p^\prime)=\pi(p)+\pi(p^\prime)$ for any $p,p^\prime$. Since
geodesics in $\mathbb{R}^k$ are lines (or their intervals), the
statement follows from Corollary \ref{corosubmersion} and a direct
computation.
\endproof
Corollary \ref{corocarnotgeodesic} can also be verified by a direct
computation using Lemma \ref{carnotconnection} and
\eqref{derivative}.
\begin{remark}For Carnot groups $G$, due to a homogeneous structure in $G$
Folland and Stein (\cite[Lemma 1.40]{FS}) proved the first statement
\eqref{control} of Corollary \ref{corcontrol} holds globally with a
universal constant $C_1$.
\end{remark}
In any sub-Riemannian manifold most sub-Riemannian geodesics are not
nonholonomic geodesics and broken geodesics. It is natural to ask
whether or not a sub-Riemannian manifold $(M,\Sigma,g_c)$ admits a
complement such that its corresponding $d_{\mathcal{H}}$ satisfies
$d_c=d_{\mathcal{H}}$. So far we can give a positive answer only for
the Heisenberg group $\mathbb{H}^n$. First we have
\begin{lemma}\label{lemcarnotgeo}
In a Carnot group $G$ any nonholonomic geodesic $\gamma$ is a
shortest (sub-)Riemannian geodesic. So if $p,q\in\gamma$,
$d_c(p,q)=d_\mathcal{H}(p,q)$.
\end{lemma}
\proof By Corollary \ref{corosubmersion} and
\ref{corocarnotgeodesic}, $\gamma$ is a Riemannian geodesic in $G$.
Since the interval or the line $\pi(\gamma)$ is obviously shortest
in Euclidean space $\mathbb{R}^k$, the horizontal lift $\gamma$ must
be a shortest Riemannian geodesic in $G$. Now the statement follows
from $d_r\leq d_c\leq d_\mathcal{H}$ where $d_r$ is the Riemannian
distance.

\endproof
\begin{proposition}
In the Heisenberg group $\mathbb{H}^n$ we have
$d_c=d_{\mathcal{H}}$.
\end{proposition}
\proof It is enough to prove that for any $\lambda>0$ and any
$p^1,p^2\in\mathbb{H}^n$ and any shortest sub-Riemannian geodesic
$\gamma$ connecting $p^1,p^2$ there exists a broken geodesic $\beta$
connecting $p^1,p^2$ such that the difference between their lengths
is less than $\lambda$.

For $p^i=(x^i,y^i,t^i)$, let
$(x^i,y^i)=\bar{p}^i=\pi(p^i)\in\mathbb{R}^{2n},i=1,2$ and
$\bar{\gamma}=\pi(\gamma)$ where $\pi$ as in the proof of Lemma
\ref{carnotsubmersion}. If $\bar{\gamma}$ is an interval of a line,
then $\gamma$ is a nonholonomic geodesic. By Lemma
\ref{lemcarnotgeo} $d_c(p^1,p^2)=d_\mathcal{H}(p^1,p^2)$.

If $\bar{\gamma}$ is an ``arc" of a ``circle", let
$A(\bar{\gamma})=\int_{\bar{\gamma}}(ydx-xdy)$ be the symplectic
area enclosing by $\bar{\gamma}$. Given $\epsilon>0$ we always can
find in $\mathbb{R}^{2n}$ a broken geodesic $\bar{\alpha}$
connecting $\bar{p}^1$ and $\bar{p}^2$ such that
$|A(\bar{\gamma})-A(\bar{\alpha})|<C\epsilon$ and
$|\ell(\bar{\alpha})-\ell(\bar{\gamma})|<\epsilon$ where $C$ is a
constant depending only on $p^1$ and $p^2$. Let $\alpha$ be the
unique horizontal lift of $\bar{\alpha}$ through $p^1$. Then
$\alpha$ is a broken geodesic in $\mathbb{H}^n$. Denote by
$p^\prime=(x^2,y^2,t^\prime)$ the other endpoint of $\alpha$. Then
from \eqref{horizontal} we get
$|t^\prime-t^2|=|A(\bar{\gamma})-A(\bar{\alpha})|<C\epsilon$. So
$d_c(p^\prime,p^2)=2\sqrt{2\pi}\sqrt{|t^\prime-t^2|}<C^\prime\sqrt{\epsilon}$.
Now we use Corollary \ref{corcontrol} or \cite[Lemma 1.40]{FS} to
get a broken geodesic $\delta$ connecting $p^\prime$ and $p^2$ with
$\ell(\delta)<C^3\sqrt{\epsilon}$ for a  constant $C^3$ depending
only on $p^1$ and $p^2$. Concatenate $\alpha$ and $\delta$ to get a
broken geodesic $\beta$ connecting $p^1$ and $p^2$. The difference
between the lengths of $\gamma$ and $\beta$ satisfies
\begin{align*}
|\ell(\beta)-\ell(\gamma)|&=|\ell(\alpha)+\ell(\delta)-\ell(\bar{\gamma})|\\
                          &=|\ell(\bar{\alpha})-\ell(\bar{\gamma})+\ell(\delta)|\\
                          &\leq \epsilon+C^3\sqrt{\epsilon}.
\end{align*}
$\beta$ is the desired broken geodesic.
\endproof

\section{Geodesically convex functions on sub-Riemannian manifolds}
The theory of convex functions on Riemannian manifolds plays a very
important role in the study of analysis and geometry on manifolds.
We expect a similar theory in the setting of sub-Riemannian setting
and wish that it could help the study of folliation theory, contact
geometry, Riemannian submersions and so on. A function on a
Riemannian manifold is convex if so is its restriction to any
Riemannian geodesic. In a similar way we define convex functions on
sub-Riemannian manifolds.

\begin{definition}[n-convex sets] A subset $\Omega$ of a sub-Riemannian
manifold $(M,\Sigma,g_c)$ is called nonholonomically geodesic convex
(n-convex for short) if for any two points $p,q\in \Omega$ and if
there exists a nonholonomic geodesic connecting them, then the
segment between them is also in $\Omega$.
\end{definition}
\begin{definition}[n-convex functions] Let $\Omega$ be a n-convex subset of a
sub-Riemannian manifold $(M,\Sigma,g_c)$. A function defined on
$\Omega$ is called n-convex if its restriction to any nonholonomic
geodesic contained in $\Omega$ is convex.
\end{definition}
 From Corollary \ref{corocarnotgeodesic} we have
\begin{proposition}In Carnot groups n-convex functions are the same
as h-convex functions defined by Danniell-Garofalo-Nieuhn
\cite{DGN}.
\end{proposition}
 By Theorem \ref{thmrieman} every Riemannian convex function on a
sub-Riemannian manifold satisfying the condition of Theorem
\ref{thmrieman} is n-convex.

\begin{theorem}A smooth function f on a sub-Riemannian manifold
$(M, \Sigma,g_c)$ is n-convex if and only if its horizontal Hessian
is positive semidefinite.
\end{theorem}
\proof Since $\mathrm{Hess}^\mathcal{H} f$ is a tensor, the
nonnegativity of $\mathrm{Hess}^\mathcal{H} f$ is equivalent to that
of $\mathrm{Hess}^\mathcal{H} f(\dot{\gamma},\dot{\gamma})$ for any
nonholonomic geodesic $\gamma$. Let $\gamma$ be a nonholonomic
geodesic. Because
\begin{align*}
\frac{d^2}{dt^2}(f\circ\gamma)&=\dot{\gamma}(\dot{\gamma}(f))(\gamma)\\
                              &=(\nabla_{\dot{\gamma}}(df))(\dot{\gamma})\\
                              &=\mathrm{Hess}
f(\dot{\gamma},\dot{\gamma})-(\nabla_{\dot{\gamma}}\dot{\gamma})f
\end{align*}
and by \eqref{hessianrep}
\begin{align*}
 \mathrm{Hess}^\mathcal{H} f(\dot{\gamma},\dot{\gamma})
 &=\mathrm{Hess}
    f(\dot{\gamma},\dot{\gamma})-B(\dot{\gamma},\dot{\gamma})\\
 &=\mathrm{Hess}
 f(\dot{\gamma},\dot{\gamma})-(\nabla_{\dot{\gamma}}\dot{\gamma})f+(D_{\dot{\gamma}}\dot{\gamma})f,
\end{align*}
 the statement follows from
$D_{\dot{\gamma}}\dot{\gamma}=0$ and the fact that
$\frac{d^2}{dt^2}(f\circ\gamma)(t)=0$ if and only $f\circ\gamma$ is
convex.
\endproof
\begin{definition}A sub-Riemannian manifold $(M,\Sigma,g_c)$
is nonholonomically complete if every nonholonomic geodesic can be
defined on the whole real line.
\end{definition}
Given a Riemannian submersion $\pi:(M,g)\rightarrow (B,g^\prime)$,
if $(B,g^\prime)$ is complete, then $(M,\Sigma,g|_{\Sigma})$ is
nonholonomically complete.
\begin{proposition}\label{proconstant}Let $(M,\Sigma,g_c)$ be a nonholonomically complete sub-Riemannian
manifold with $\Sigma$ regular and $M$ connected. If a n-convex
function $f:M\rightarrow\mathbb{R}$ is upper bounded, then $f$ must
be a constant.
\end{proposition}
\proof Let $f:M\rightarrow\mathbb{R}$ be a n-convex function. Since
$M$ is nonholonomically complete, every nonholonomic geodesic
$\gamma(t)$ can be defined in $(-\infty,\infty)$. Let $p,q$ be in a
nonholonomic geodesic $\beta(t):[0,1]\rightarrow M$ such that
$\beta(0)=p$ and $\beta(1)=q$. Extend $\beta$ to $[0,t]$ for $t>1$.
By setting $u=st,s\in [0,1]$, $\beta$ is reparameterized to
$\bar{\beta}(s):=\beta(st)$. The convexity of $f$ on $\bar{\beta}$
implies
\begin{align*}
f(\bar{\beta}(s))&\leq(1-s)f(\bar{\beta}(0))+sf(\bar{\beta}(1))\\
                 &=(1-s)f(p)+sf(\beta(t))
\end{align*}
for any $s\in[0,1]$. In particular for $s=\frac{1}{t}$ we get
\begin{align*}
 f(q)&=f(\beta(1))=f(\bar{\beta}(\frac{1}{t}))\\
      &\leq(1-\frac{1}{t})f(p)+\frac{1}{t}f(\beta(t)).
\end{align*}
Since $f$ is upper bounded, letting $t$ go to $\infty$ we obtain
$f(p)\leq f(q)$. Similarly we have $f(p)\geq f(q)$. Thus $f$ is
constant on any nonholonomic geodesic. Now the statement follows
from Theorem \ref{connectivity}.
\endproof

We turn to the regularity problem of n-convex functions on regular,
connected sub-Riemannian manifolds. Balogh and Rickly (\cite{BR})
proved that for the Heisenberg group $\mathbb{H}^n$ h-convex
functions are locally Lipschitz with respect to the
Carnot-Carath\`{e}odory distance. The arguments and regularity
result in \cite{BR} were later extended to Carnot groups of step 2
by \cite{Ri} and \cite{SY} independently. For general Carnot groups,
\cite{Ma} proved that any h-convex function with local upper bound
is locally Lipschitz. For more properties, in particular the
equivalence of several notions of h-convex functions, and their
applications in nonlinear subelliptic  PDEs, we refer to
\cite{DGN,LMS,DGNT,GT,GM,JLMS,BR,Ma,Wang,Ri} and references therein.
In the following we assume $\Omega$ be a n-convex open subset in a
sub-Riemannian manifold $(M,\Sigma,g_c)$ with $\Sigma$ regular and
$M$ connected.
\begin{proposition}\label{bound}Any n-convex function
$f:\Omega\rightarrow\mathbb{R}$ locally bounded above is also
locally bounded below. More explicitly, if $f\leq C$ in a
neighborhood of $y_0$, then there exists a constant
$\bar{\epsilon}_0$ such that $f(p)\geq 2^Nf(y_0)-(2^N-1)C$ for any
$p\in{B(y_0,\bar{\epsilon}_0)}$, where $N$ is the integer in
Corollary \ref{corcontrol}.
\end{proposition}
\proof For $y_0\in{\Omega}$ let $\mathcal{O}$ be the neighborhood of
$y_0$ as in Theorem \ref{frame} where there exists a horizontal
frame $\{E_1,\cdots,E_k\}$ of norm 1 such that integral curves of
$E_i$ are nonholonomic curves. Let $\mathcal{O}^\prime$ be a
neighborhood of $y_0$ such that
$\bar{\mathcal{O}}^\prime\subset\mathcal{O}$ and there exists
$\tilde{\epsilon}_0>0$ such that for any $p\in{\mathcal{O}^\prime}$
and $i=1,\cdots,k$, the integral curve $\Psi_i^p(t)$ of $E_i$ with
$\Psi_i^p(0)=p$ is defined in the interval
$[-\tilde{\epsilon}_0,\tilde{\epsilon}_0]$. Take
$\tilde{\epsilon}_0$ smaller if necessary such that
$B(y_0,\tilde{\epsilon}_0)\subset\mathcal{O}^\prime$. Let
$B(y_0,\bar{\epsilon}_0)\subset\mathcal{O}^\prime$ be the ball as in
Corollary \ref{corcontrol}. Again take $\bar{\epsilon}_0$ smaller if
necessary such that $NC_1\bar{\epsilon}_0<\tilde{\epsilon}_0$.

By Corollary \ref{corcontrol} for any
$p\in{B(y_0,\bar{\epsilon}_0)}$, there exists a broken geodesic
$\gamma$ connecting $y_0$ and $p$, which consists of $N$
nonholonomic geodesics, each one of length less than
$C_1d_c(y_0,p)$. Let $q_i (i=1,\cdots,N-1)$ be the $N-1$ breaks and
$\gamma_0$ the nonholonomic geodesic from $y_0$ to $q_1$, $\gamma_i$
the one from $q_{i-1}$ to $q_i$ and finally $\gamma_{N-1}$ from
$q_{N-1}$ to $p$. From our choices of $\bar{\epsilon}_0$ and
$\tilde{\epsilon}_0$, we have $\gamma_i\subset\mathcal{O}^\prime$
and $\ell(\gamma_i)<\tilde{\epsilon}_0$, $i=1,\cdots,N-1$. Because
$E_i$'s are of norm 1, $\gamma_i$'s are parameterized by length and
thus $\gamma_i$ can be extended to $\bar{\gamma}_i$ such that $q_i$
is the middle point of $\bar{\gamma}_i$, $i=0,\cdots,N-1$ with
$q_0=y_0$. Also $\bar{\gamma}_i\subset\mathcal{O}^\prime$.

Assume the n-convex function $f$ in $\mathcal{O}^\prime$ is bounded
by $C$ from above. We claim $f$ is bounded below in
$B(y_0,\bar{\epsilon}_0)$. In fact, because by definition $f$ is
convex on each $\bar{\gamma}_i$ and $q_i$ is the middle point of
$\bar{\gamma}_i$ $i=0,\cdots,N-1$, we have for $p\in
B(y_0,\bar{\epsilon}_0)$
\begin{align*}
f(p)&\geq2f(q_{N-1})-C\\
    &\geq2(2f(q_{N-2})-C)-C=2^2f(q_{N-2})-(1+2)C\\
    &\cdots\\
    &\geq2^{N-1}f(q_1)-(1+2+2^2+\cdots+2^{N-2})C\\
    &\geq2^{N-1}(2f(q_0)-C)-(1+2+2^2+\cdots+2^{N-2})C\\
    &=2^Nf(y_0)-(1+2+2^2+\cdots+2^{N-1})C\\
    &=2^Nf(y_0)-(2^N-1)C
\end{align*}
\endproof
Now we prove the locally Lipschitz continuity of n-convex functions
locally bounded above. We follow closely the arguments in \cite{Ri}
for the case of Carnot groups.
\begin{theorem}\label{thmmain}If a n-convex function $f:\Omega\rightarrow\mathbb{R}$
is locally bounded above, then $f$ is locally Lipschitz continuous
with respect to $d_c$ and $d_\mathcal{H}$.
\end{theorem}
\proof Let
$y_0,\mathcal{O},E_i,\mathcal{O}^\prime,\tilde{\epsilon}_0,\bar{\epsilon}_0$
be as in the proof of Proposition \ref{bound} such that
$(2NC_1+1)\bar{\epsilon}_0<\frac{1}{2}\tilde{\epsilon}_0$. Let $f$
be bounded in $B(y_0,\tilde{\epsilon}_0)$, $|f|\leq C$. We first
prove that $f$ is Lipschitz on any integral curve $\gamma$ of
$E_i$'s which is contained in
$B(y_0,\frac{1}{4}\tilde{\epsilon}_0)$. To this aim we assume
$\gamma(0)\in B(y_0,\frac{1}{4}\tilde{\epsilon}_0)$. Since $\gamma$
is parameterized by length and is contained in
$B(y_0,\frac{1}{4}\tilde{\epsilon}_0)$,
$\ell(\gamma)\leq\frac{1}{2}\tilde{\epsilon}_0$. According to the
choice of $\tilde{\epsilon}_0$, $\gamma$ can be extended. Still
denote by $\gamma(t)$ the extended $\gamma$. Set
$t_-:=\max\{t<0|\gamma(t)\in\partial
B(y_0,\frac{1}{2}\tilde{\epsilon}_0)\}>-\tilde{\epsilon}_0$ and
$t_+:=\min\{t>0|\gamma(t)\in\partial
B(y_0,\frac{1}{2}\tilde{\epsilon}_0)\}<\tilde{\epsilon}_0$. Then
$t_+-t_-\in [\frac{1}{2}\tilde{\epsilon}_0,\tilde{\epsilon}_0]$, and
if $t\in[t_-,t_+]$ and $\gamma(t)\in
B(y_0,\frac{1}{4}\tilde{\epsilon}_0)$,
$t-t_-\geq\frac{1}{4}\tilde{\epsilon}_0$ and
$t_+-t\geq\frac{1}{4}\tilde{\epsilon}_0$. Thus we can pick
$\lambda\in{[\frac{1}{4},\frac{3}{4}]}$ such that
$t=(1-\lambda)t_-+\lambda t_+$. Let $t_1,t_2\in{[t_-,t_+]}$ such
that $t_1<t_2$ and
$\gamma(t_1),\gamma(t_2)\in{B(y_0,\frac{1}{4}\tilde{\epsilon}_0)}$.
Then $t_i=(1-\lambda_i)t_-+\lambda_i t_+$ where
$\lambda_i\in{[\frac{1}{4},\frac{3}{4}]},i=1,2 $ and
$\lambda_1<\lambda_2$. So
$$
t_1=\frac{\lambda_2-\lambda_1}{\lambda_2}t_-+\frac{\lambda_1}{\lambda_2}t_2\quad\textrm{
and }\quad
t_2=\frac{1-\lambda_2}{1-\lambda_1}t_1+\frac{\lambda_2-\lambda_1}{1-\lambda_1}t_+.
$$
Because as an integral curve of some $E_i$, $\gamma$  is a
nonholonomic geodesic, from the convexity of $f$ we obtain
\begin{align*}
f(\gamma(t_1))-f(\gamma(t_2))&\leq\frac{\lambda_2-\lambda_1}{\lambda_2}f(\gamma(t_-))+
                             \frac{\lambda_1-\lambda_2}{\lambda_2}f(\gamma(t_2))\\
                             &=\frac{t_2-t_1}{\lambda_2(t_+-t_-)}\{f(\lambda(t_-)-f(\gamma(t_2)))\}\\
                             &=\frac{d_\mathcal{H}(\gamma(t_1),\gamma(t_2))}{\lambda_2(t_+-t_-)}\{f(\gamma(t_-))
                               -f(\gamma(t_2))\}\\
                             &\leq\frac{8C}{\tilde{\epsilon}_0}d_c(\gamma(t_1),\gamma(t_2)),
\end{align*}
where we use $d_\mathcal{H}\leq d_c$ and the fact that when
$\tilde{\epsilon}_0$ is small enough, $\gamma$ is the unique
nonholonomic geodesic passing through $\gamma(0)$. Similarly we have
$$
f(\gamma(t_2))-f(\gamma(t_1))\leq
\frac{8C}{\tilde{\epsilon}_0}d_c(\gamma(t_1),\gamma(t_2)).
$$
We have proved that
$$
|f(p_1)-f(p_2)|\leq\frac{8C}{\tilde{\epsilon}_0}d_c(p_1,p_2)
$$
for any $p_1,p_2\in{\gamma}\subset
B(y_0,\frac{1}{4}\tilde{\epsilon}_0)$ where $\gamma$ is a segment of
an integral curve of some $E_i$.

Now for any $p,q\in{B(y_0,\bar{\epsilon}_0)}$, by Corollary
\ref{corcontrol} there exists a broken geodesic $\beta$ connecting
$p$ to $q$, and $\beta$ consists of $N$ nonholonomic geodesics, as
an integral curve of some $E_i$ each one of length less than
$C_1d_c(p,q)$. Since
\begin{align*}
d_c(y_0,\beta)&\leq d_c(\beta,p)+d_c(y_0,p)\\
               &\leq NC_1d_c(p,q)+d_c(y_0,p)
               \leq(2NC_1+1)\bar{\epsilon}_0\\
               &\leq \frac{1}{2}\tilde{\epsilon}_0,
\end{align*}
each nonholonomic geodesic segement of $\beta$ is contained in
$B(y_0,\tilde{\epsilon}_0)$. Thus to each segment  we can use the
above estimate to get
$$
|f(p)-f(q)|\leq \frac{8NCC_1}{\tilde{\epsilon}_0}d_c(p,q)\leq
\frac{8NCC_1}{\tilde{\epsilon}_0}d_\mathcal{H}(p,q).
$$
\endproof
From Theorem \ref{thmmain} and Proposition \ref{proconstant} or
\ref{bound} we have
\begin{corollary} Let $(M,\Sigma,g_c)$ be a compact sub-Riemannian
manifold with $\Sigma$ regular and connected. Then on $M$ any
n-convex function locally bounded above is constant.
\end{corollary}

\end{document}